\documentclass[a4paper,12pt, legno]{article}
\usepackage{amsmath}
\usepackage{amsfonts}
\usepackage{amsthm}
\usepackage{mathrsfs}
\usepackage[T1]{fontenc}
\usepackage{amsmath}
\usepackage{amsfonts}
\usepackage{amsthm}
\usepackage{mathrsfs}
\usepackage{amssymb}
\usepackage{pst-node}
\usepackage{tikz-cd}

\newtheorem{prop}[]{Proposition}

\newtheorem{cor}{Corollary}

\newtheorem{rem}[]{Remark}

\newcommand{\Rset}{\mathbb{R}}
\newcommand{\Cset}{\mathbb{C}}

\begin{document}

\title{On hyperbolic characteristic functions from an analytic  and a free-probability point of view}
\author{Zbigniew J.  Jurek (University of Wroc\l aw)}
\date{May 25, 2020}
\maketitle

\begin{quote} \textbf{Abstract.} For  free-probability Voiculescu  transforms, analogous to   hyperbolic characteristic functions,  we show how to get their representing measures in an  integral form. For that purpose it is enough to know those transforms  only  on the  imaginary  axis. This  is in a  contrast to a complex analysis where one needs to know them in some domains in the complex plane.

\medskip



\medskip
\emph{Abbreviated title: On hyperbolic characteristic function.}

\end{quote}

\underline{Author's address}: 

Institute of Mathematics, University of Wroc\l aw, Pl. Grunwaldzki 2/4, 

50-384 Wroc\l aw, Poland;

Email:  zjjurek@math.uni.wroc.pl ;   www.math.uni.wroc.pl/$\sim$zjjurek

\medskip
\medskip
In a  classical  complex analysis one of  the  fundamental results is the  integral representation of \emph{ analytic functions}
from the upper to the  lower complex half-plane. Those functions, say  $H$, admit  unique canonical form
$$
H(z)= a + \int_{\Rset} \frac{1+z x}{z-x}\rho(dx)= a -\int_{\Rset}\big[ \frac{1}{x-z}- \frac{x}{1+x^2} \big] (1+x^2)\rho(dx) ,  \  (\star)
$$
where $a\in \Rset$ is a constant, $\rho$ is a finite (Borel) measure on a real line $\Rset$; (in a literature $H$ are  coined  as  Pick functions and a representation $(\star)$ is called  Nevalinna Theorem.)  
One simply notes that  a constant $a=\Re H(i)$ (a real part), and a   total mass $\rho(\Rset)= - \Im (H(i)$.  Finally, for a measure $\rho$ we have an inversion formula
$$
\rho([c,d])= \lim_{\epsilon \to 0^+}\frac{1}{\pi}\int_c^d \Im H(x+ i \epsilon)dx, \ \mbox{whenever } \  \rho(\{c,d\})=0;  \qquad \ (\star \star)
$$
cf.Akhiezer(1965), p.125, or Lang(1975), p.380, or Bondensson(1992), p.21. 

\medskip
However, what can be said  about a measure  $\rho$  if we only  have values $H(it)$, for $t \neq 0,$   and  we  don't know if  it is  a restriction of  an analytic function to the  imaginary axis ? Note that in $(\star \star)$ we need to know $H$ in some strips  of the complex plane  to retrieve  a measure $\rho$. Never the less,  in Jankowski and Jurek (2012), Theorem 1, there is an inversion procedure
that allows us to identify a measure $\rho$ or more precisely its characteristic function $\hat{\rho}$.

In a couple of last decades representation of the form $(\star)$ appeared in so-called \emph{free-probability}  as a free-analog of the classical  L\'evy-Khintchine formula for infinite divisible characteristic functions (Fourier transforms). 

In this  note we show applications of the  inversion formula from Jankowski and Jurek (2012) for  $\tilde{C},\tilde{S}, \tilde{T}$  free-analogs of the classical hyperbolic functions $C,S$ and $T$. Recall that $C, S$ and $T$  are defined by  their characteristic functions as follows
$$\phi_C(t):=\frac{1}{\cosh(t)},\   \   \phi_S(t):=\frac{\sinh(t)}{t}, \  \ \phi_T(t):= \frac{\tanh(t)}{t}, \ t \in \Rset.$$   
In free-probability variables $\tilde{C},\tilde{S}, \tilde{T}$ are given by their Voiculescu transforms $V_{\tilde{C}}(z), V_{\tilde{S}}(z)$ and $ V_{\tilde{T}}(z), z \in \Cset^+$, although in our approach to free-probability theory we consequently use only purely imaginary $z=it, t\neq 0$; cf. Jurek (2019), Corollaries 3 ,4 and 5, respectively.

\medskip
Let us recall that hyperbolic characteristic functions, from an infinite divisibility point of view,  were studied  in Pitman and Yor (2003) and from a  selfdecomposability point of view in  Jurek (1997) (as infinite series of independent exponentially distributed variables), and in Jurek-Yor (2004) (from stochastic representations of their background driving processes). The last description   can be done  as all hyperbolic characteristic functions are selfdecomposable ones  and therefore they  admit  a representation  by random integrals; 
cf. Jurek and Mason (1993), Chapter 3, Theorem 3.6.8 or  Jurek and Vervaat (1983), Theorem 3.2.

\medskip
\medskip
\textbf{1. Introduction.} 

\medskip
Let us for an index $X$  (where $X$ can be a random variable or a measure or a characteristic function) define a function  $V_X$ on  the  imaginary  axis   $ i(\Rset\setminus{\{0\}})$  as follows
 \begin{equation}
V_X(it):= a_X +\int_{\Rset}\frac{1+itx}{it-x} m_X(dx), \ \ t\neq 0,
\end{equation}
where $a_X\in \Rset$  and  $m_X$ is a non-negative, finite Borel measure. Then note that   $V_X(i)=a_X -  i\,m_X(\Rset)$ and hence we get  
\begin{equation}
  a_X=\Re V_X(i) \in \Rset , \  \  \   m_X(\Rset)= - \Im V_X(i) \in [0,\infty).
\end{equation}
Furthermore, if  $\hat{m}_X(s):= \int_\Rset e^{isx} m_X(dx)$ denotes \emph{a characteristic function of a measure $m_X$} then its \emph{Laplace transform $\mathfrak{L}[ \hat{m}_X; w]$, for $w>0$}, satisfies  equality
\begin{equation}
\mathfrak{L}[\hat{m}_X; w] = \int_0^\infty\hat{m}_X(x)e^{- w x}dx = \frac{i V_X(-i w) - i \Re V_X(i) - w \Im V_X(i)}{w^2-1};
\end{equation}
cf. Theorem 1 in Jankowski  and Jurek (2012).  
Equivalently,
\begin{multline}
\mathfrak{L}[\hat{m_X}(s)+i a_X \sinh(s) - m_X(\Rset)\cosh (s) ; w] =  \frac{ iV_X(-iw)}{w^2-1}, \ \ \ 
\mbox{because} \\
\mathfrak{L}\big[\sinh x; w]= \frac{1}{w^2-1}, \  \mbox{and}  \   \mathfrak{L}[\cosh x; w]=\frac{w}{w^2-1},  \ \mbox{for}  \  w>1.
\end{multline}

\medskip
In propositions below our aim is to show that  functions  $ w\to \frac{i V_X(-i w)}{w^2-1}$,

\noindent  (on the right hand side in (4)),  are indeed  Laplace's transform of some  functions or  
measures. This , in principle, enables us to identify a representing measure $m_X$ in the canonical form (1).

\medskip
\medskip
\textbf{2. Results.}
In order to present our results we will need some special functions.. Therefore before each proposition we recall   their definitions and/or  characterizations. Many of  those  functions 
are derived from \emph{Euler's $\Gamma$ gamma function}: $\Gamma(z):=\int_0^\infty x^{z-1}e^{-x}dx; \ \Re z>0$ and \emph{digamma  function} $\psi(z):= d/dz\log\Gamma(z)$.
 For more  facts and  formulas we refer to the Appendix at the end of this article.

\medskip
In the first proposition we need  \emph{$\beta$ beta function} which admits a representation  $\beta(z)=\int_0^\infty (1+e^{-x})^{-1} \, e^{-z x} dx, \  \ \Re z >0$, and originally was
defined via digamma function $\psi$; cf. Appendix (A).

\begin{prop}  
For a free-infinitely divisible Voiculescu transform $$ V_{\tilde{C}}(it)=i[1-t \beta(t/2)], t\neq 0,$$  we have that  in its representation (1) a real parameter $a_{\tilde{C}}=0$ 
and a measure $m_{\tilde{C}}$ is  such that a total mass $m_{\tilde{C}}(\Rset)= \pi/2-1\approx 0.57079$  and its  characteristic  function $\hat{m}_{\tilde{C}}$  is equal to 
\begin{multline}
\hat{m}_{\tilde{C}}(s)=2\sinh(s)\tan^{-1}(e^{-s})+\frac{\pi}{2} e^{-s}-1 \\  =  
\int_0^\infty\cos(sx)\frac{|x|}{1+x^2}\frac{1}{\sinh(\pi|x|/2)} dx
\end{multline}
($\tilde{C}$ indicates a free-probability analogue of a classical hyperbolic cosine characteristic function $\phi_C(t)=1/\cosh(t).$)\end{prop}

To formulate next proposition we need  two special functions. Namely, \emph{ a digamma function} $\psi$ , which is defined as $\psi(z):= \frac{d}{dz}\log \Gamma(z)$,  and \emph{an exponential integral function}
$Ei(x):=-\int_{-x}^\infty\frac{e^{-t}}{t}dt$, for $ x <0$; for more see Appendix (B).

\begin{prop}.
For a free-infinitely divisible Voiculescu transform 
$$V_{\tilde{S}}(it)= i[t\psi(t/2)-t\log(t/2) +1] , t\neq 0,$$
we have that  in its  representation (1)  a parameter $a_{\tilde{S}}=0$ and a measure $m_{\tilde{S}}$ is  a such that 
$m_{\tilde{S}}(\Rset)=\gamma + \log2-1\approx 0.270362$  
and its characteristic function $\hat{m}_{\tilde{S}}$  is of a form
\begin{multline}
\hat{m}_{\tilde{S}}(s)  =   -1 +\cosh(s)\big( \log(1+e^{-s}) - \log(1-e^{-s})\big)+\frac{e^{-s}}{2}Ei(s)+ \frac{e^{s}}{2} Ei(-s) \\ 
 = 2 \int_0^\infty \cos(sx)\frac{x}{1+x^2}\frac{1}{e^{\pi x}-1}dx, \  \ s>0, \qquad \qquad 
\end{multline}
and $Ei(x)$ is an exponential integral function.

($\tilde{S}$ indicates a free-probability counter part of the hyperbolic sine characteristic function $\phi_S(t)= t/ \sinh(t)$.)
\end{prop}

\medskip
In next statements appear special functions from  two previous propositions because of the  elementary relation: $\phi_C(t)=\phi_S(t)\cdot \phi_T(t)$.
\begin{prop}
For a free-infinitely divisible Voiculescu transform 
$$V_{\tilde{T}}(it)= V_{\tilde{C}}(it)- V_{\tilde{S}}(it) = it\,[\,\log (t/2)-\beta(t/2)-\psi(t/2)\,], t\neq 0,$$
we have that  in its  a representation (1)  a parameter $a_{\tilde{T}}=0$ and for a measure $m_{\tilde{T}}$ we have that $m_{\tilde{T}}(\Rset)= \pi/2- \gamma -\log2\, \, (\approx 0.3004)$, and its characteristic function $\hat{m}_{\tilde{T}}$  has a form
\begin{multline}
\hat{m}_{\tilde{T}}(s)=   \frac{\pi}{2} e^{-s}+ 2 \sinh(s)\tan^{-1}(e^{-s}) - \cosh s\log(1+e^{-s})  \\  -  \frac{e^{-s}}{2}\big(Ei(s)-\log(1-e^{-s}) \big) - \frac{e^{s}}{2}\big(Ei(-s)-\log(1-e^{-s}) \big) 
\\ =\int_0^\infty \cos(sx)\frac{|x|}{1+x^2}\frac{e^{-\pi|x|/4}}{\cosh (\pi|x|/4)}dx,   \ \  s>0.  \qquad \qquad 
\end{multline}
($T$ indicates a hyperbolic  tangent  characteristic function $\phi_T(t)= \tanh(t)/t$.)
\end{prop}

\medskip
\medskip
All three Voiculescu  transforms  from Propositions 1, 2 and 3 are  free-probability analogs of selfdecomposable characteristic functions. Therefore  they have so called \emph{background driving 
terms} from corresponding random integral representations. In particular, these are  background driving characteristic  functions $\psi_{\tilde{C}}, \psi_{\tilde{S}}$ and $\psi_{\tilde{T}}$, and L\'evy (spectral) measures $N_C,N_S$ and $N_T$; cf. Jurek (2019), Section 2.1  or  Jurek-Yor (2004).

As in previous propositions we have similar results for them as well, although we computed  it only for a background driving characteristic function  $\psi_{\tilde{C}}$, only;  cf. Proposition 4 below.

\medskip
For a following proposition,  we need another two special functions. Namely, \emph{Riemann's zeta  function}, $\zeta$,  and \emph{polylogarithm functions $Li_n(z)$ } (in Wolframalpha.com language  written as  $PolyLog[n,z]$); cf. Appendix (C).
\begin{prop}
For a free-infinitely divisible Voiculescu transform
$$
V_{\psi_{\tilde{C}}} (it) = i\,[\frac{t^2}{2}\zeta(2, t/2)-\frac{t^2}{4}\zeta(2, t/4)+1]
$$
we have that in its representation (1) a parameter $a_{\psi_{\tilde{C}}}=0$, a total mass $m_{\psi_{\tilde{C}}}(\Rset)= 2C-1\,   (\approx 0.83193)$ and for a measure $m_{\psi_{\tilde{C}}}$ we have its characteristic function
\begin{multline*}
\hat{m}_{\psi_{\tilde{C}}} (t)= -1 - t \tanh(t)  - \cosh(t) \big ( \,i(Li_2(ie^t) - Li_2(-i e^t))  +2 t \arctan(e^t)\,\big) \\ =
\frac{\pi}{2} \int_0^\infty \cos(tx)\,\frac{x^2}{1+x^2}\,\frac{\cosh(\pi x/2)}{\sinh^2(\pi x/2)}dx. \qquad \qquad 
\end{multline*}
In particular, we have:
$$   i ( Li_2(i e^t)- Li_2(-i e^t)) = 2\sum_{k=1}^\infty (-1)^k \frac{e^{(2k-1)t}}{(2k-1)^2}; \  \  i ( Li_2(i )-Li_2(-i))= - 2C. $$
\end{prop}
As a by-product of our Proposions 1 and 4, we have
\begin{cor}
For a hyperbolic cosine function $\phi_C(t)= 1/\cosh t$ we have
\[
\hat{m}_{\psi_{\tilde{C}}}(t) +\hat{m}_{\tilde{C}}(t)=2\int_0^\infty\cos(tx)\frac{x^3}{1+x^2}\big( k_C(x)\big)^\prime dx
\]
\noindent where a function $k_C(x):= ( 2x \sinh (\pi x/2))^{-1}$ is a denstity of L\'evy measure of  a hyperbolic cosine  function. $\phi_C$.
\end{cor}
\medskip
\textbf{3. PROOFS}

\medskip
All boldface numbers appearing below refer to the corresponding fromulae in  Gradstheyn-Ryzhik (1994).

\medskip
\underline{\emph{Proof of Proposition 1.} }

\medskip
First,  note that $V_{\tilde{C}}(i)=i (1-\beta(1/2))= -i(\pi/2-1)$ and therefore  $m_{\tilde{C}}(\Rset)=\pi/2-1$. Second, since by (8),
\[
\beta(s)= \int_0^\infty \frac{1}{1+e^{-x}}e^{-s x}dx=\mathfrak{L}[ (1+e^{-x})^{-1}; s], \ \ \Re s>0, \  \textbf{8.371}(2)
\]
consequently
\begin{multline}
\frac{iV_{\tilde{C}}(-iw)}{w^2-1}=\frac{1- w\beta(w/2)}{w^2-1}=\mathfrak{L}[\sinh x; w] - \mathfrak{L}[\cosh x;w]\,\mathfrak{L}[2/(1+e^{-2x});w]\\ =\mathfrak{L}[\sinh x; w] -\mathfrak{L}[\,\underline{(\cosh(s)\ast (2/(1+e^{-2s})))}(x); w], \ \ \qquad
\end{multline}
where $\ast$ denotes a convolution of functions on positive half-line.

Third, one checks by a differentiation (or by WolframAlpha or Mathematica) that for $x>0$ we have
\begin{multline*}
\big(\underline{\cosh(s)\ast (2/(1+e^{-2s}))}\big)(x):=\int_0^x\cosh(x-s)\,\frac{2}{1+e^{-2(x-s)}}ds \\ =2\sinh(x)\arctan(e^{s-x})-e^{-s} + constant |_{s=0}^{s=x} \\ = 2 \sinh(x)\arctan (e^{x-x})- e^{-x} - \big( 2\sinh(x) \arctan(e^{-x}) -1\big)\\ = 2\sinh(x)(\pi/4 - \arctan(e^{-x}))-e^{-x}+1, \qquad \qquad
\end{multline*}
and inserting it into (9) we get
\begin{equation}
\frac{iV_{\tilde{C}}(-iw)}{w^2-1} = \mathfrak{L}[\sinh x -2\sinh x (\pi/4-\arctan(e^{-x}))+e^{-x}-1;w]
\end{equation}
Finally, since $m_C(\Rset)=\pi/2-1$ and taking into account (4) we have
\begin{multline*}
\hat{m}_C(x) = (\pi/2-1)\cosh x+ \sinh x -2\sinh x (\pi/4-\arctan(e^{-x})+e^{-x}-1 \\ = 2\sinh x \arctan(e^{-x})+\pi/2 e^{-x}-1, \qquad
\end{multline*}
which gives a first equality in (5).

On the other hand,
from Jurek (2019), Example 1 we know that $V_{\tilde{C}}(it)$ is a free-probability analog of a classical hyperbolic characteristic function $1/\cosh(t)$ whose (finite) Khintchine measure $m_C$ has a density 
$\frac{1}{2}\frac{|x|}{1+x^2}\frac{1}{\sinh(\pi|x|/2)}$, for $x\in \Rset$. Thus
\[
\hat{m}_C(t)=\int_{\Rset} e^{itx}\frac{1}{2}\frac{|x|}{1+x^2}\frac{1}{\sinh(\pi|x|/2)}dx= \int_0^\infty \cos (tx)\frac{|x|}{1+x^2}\frac{1}{\sinh(\pi|x|/2)}dx,
\]
which completes  proof of Proposition 1.
\begin{rem}
\emph{ From a formula  in Proposition 1, we have an identity 
\[
\int_0^\infty \cos ( sx)\, \frac{|x|}{1+x^2}\frac{1}{\sinh (\pi|x|/2)} dx = -1 +\frac{\pi}{2} e^{-s}+ 2 \sinh(s)\tan^{-1}(e^{-s}), \ s>0;
\]
which is confirmed by \textbf{4.113}(8).
}
\end{rem}

\medskip
\underline{\emph{Proof of Proposition 2.} }

\medskip
Since $\psi(1/2)= - \gamma -2 \log2$, (\textbf{8.366}(2))  then $V_S(i)= - i(\gamma +\log 2-1)$.  Hence  a  parameter $a_{\tilde{S}}=0$  and a finite measure $m _{\tilde{S}}$ in (1)  has a finite mass   $m_{\tilde{S}}(\Rset)= \gamma + \log 2-1$.  Finally,  using an integral formula for $\psi$ function,
\[
\psi(z)=\log z +\int_0^\infty\big(\frac{1}{s}-\frac{1}{1-e^{-s}}\big) e^{-  z s}ds, \ \Re z >0; \ \ \textbf{8.361}(8),
\]
we have 
\begin{multline}
\frac{i V_S(-iw)}{w^2 -1}= \frac{1+w (\psi(w/2)-\log(w/2))}{w^2-1} \\ =  \mathfrak{L}[\sinh x; w] +\mathfrak{L}[\cosh x; w] \mathfrak{L}[ \frac{1}{x}-\frac{2}{1-e^{-2x}}; w] 
\\ =\mathfrak{L}[\sinh x +\underline{\big(\cosh s \ast (\frac{1}{s}-\frac{2}{1-e^{-2s}})\big)}(x); w],
\end{multline}
where $"\ast"$  denotes a convolution of functions.

By a direct differentiation using properties quoted before  a proof of this  proposition (or using  www.Wolframalpha.com) one checks that:
\begin{multline*}
\int \cosh(x-s) \Big[\frac{1}{s} -\frac{2}{1-e^{-2s}}\Big] ds \\
= \frac{e^{-x}}{2}Ei(s)+\frac{e^x}{2}Ei(-s) - e^{s-x}+\cosh( x)\, \log\frac{1+e^{-s}}{1-e^{-s}}+ constant.
\end{multline*}
Thus for $x>0$ we have 
\begin{multline*}
\underline{\big(\cosh s \ast (\frac{1}{s}-\frac{2}{1-e^{-2s}}) \big)}(x)
=\int_0^x \cosh(x-s) \Big[\frac{1}{s} -\frac{2}{1-e^{-2s}}\Big] ds \\
= \big[ \frac{e^{-x}}{2}Ei(s)+\frac{e^x}{2}Ei(-s) - e^{s-x}+\cosh( x)\, \log\frac{1+e^{-s}}{1-e^{-s}}\Big]\Big|^{s=x}_{s=0^+}  \\ =
\frac{e^{-x}}{2}Ei(x)+\frac{e^x}{2}Ei(-x) - e^{x-x}+\cosh( x)\, \log\frac{1+e^{-x}}{1-e^{-x}} \\  - 
\lim_{s \to 0^+} \big [ \frac{e^{-x}}{2}(Ei(s) -\log(1-e^{-s}))+ \frac{e^x}{2}(Ei(-s)-\log(1-e^{-s}))  \\ - e^{s-x}+ \cosh x \log(1+e^{-s})\big] = \frac{e^{-x}}{2}Ei(x)+\frac{e^x}{2}Ei(-x) - 1 \\ +\cosh( x)\, \log\frac{1+e^{-x}}{1-e^{-x}} -
\gamma \cosh x +e^{-x} -\cosh x \log2\\  =\frac{e^{-x}}{2}Ei(x)+\frac{e^x}{2}Ei(-x) - 1 +e^{-x} +
\cosh(x) \big ( \log\frac{1+e^{-x}}{1-e^{-x}} -\gamma-  \log 2\big).
\end{multline*}
Inserting above line into (11) and  using  (4) with $m_S(\Rset)= \gamma + \log 2-1$ we get
\begin{multline}
\hat{m}_S(x)= (\gamma +\log2-1)\cosh(x)+\sinh(x) \\ + \frac{e^{-x}}{2}Ei(x)+\frac{e^x}{2}Ei(-x) - 1 +e^{-x} +
\cosh(x) \big ( \log\frac{1+e^{-x}}{1-e^{-x}} - (\gamma+ \log 2) \big)\\=
-\cosh(x)+\sinh(x)+ \frac{e^{-x}}{2}Ei(x)+\frac{e^x}{2}Ei(-x) - 1 +e^{-x} +
\cosh(x) \log\frac{1+e^{-x}}{1-e^{-x}} \\= \frac{e^{-x}}{2}Ei(x)+\frac{e^x}{2}Ei(-x) - 1 + \cosh(x)(\log(1+e^{-x})-\log(1-e^{-x})),
\end{multline}
which proves first equality in Proposition 2.

Since from Jurek (2019), Corollary 4,  we know that $V_{\phi_S}$ is a free-analog of a classical hyperbolic sine characteristic function $\phi_S(t)= t/\sinh(t)$ whose (finite) Khintchine measure is equal to 
\[
m_S(dx)= \frac{1}{2}\frac{|x|}{1+x^2} \frac{e^{-\pi|x|/2}}{\sinh(\pi|x|/2)}dx= \frac{|x|}{1+x^2}\, \frac{1}{e^{\pi |x|}-1}dx, \  \mbox{on} \ \Rset,
\]
we get second equality in Proposition 2.
\begin{rem}
\emph{From Proposition 2 we get  an  integral identity
\begin{multline*}
 2 \int_0^\infty \cos(sx)\frac{x}{1+x^2}\frac{1}{e^{\pi x}-1}dx \\  =   -1 +\cosh(s)\big( \log(1+e^s) - \log(1-e^{-s})\big)+\frac{e^{-s}}{2}Ei(s)+ \frac{e^{s}}{2} Ei(-s)  \  \ s>0,
\end{multline*}
that might be of some  interest and it seems to be new? }
\end{rem}

\medskip
Since for hyperbolic characteristic functions $\phi_C, \phi_S$ and $\phi_T$ we have that $\phi_C(t)= \phi_S(t)\cdot\phi_T(t)$ therefore for their Khintchine measures $m_C, m_S, m_T$ we have $m_C(dx)=m_S(dx)+m_T(dx)$.

\medskip
\underline{\emph{Proof of Proposition 3.} }

\medskip
Taking into account a discussion above and the fact that free-infinitely divisible transforms  $V_{\tilde{C}}, V_{\tilde{S}}, V_{\tilde{T}}$, in Jurek (2019)  were defined via one - to -one  correspondence with classical infinite divisibility we get that 
$V_{\tilde{C}}(it)=V_{\tilde{S}} (it)+V_{\tilde{T}}(it), t \neq 0$. Consequently,  proof of Proposition 3 follows from proofs of Propositions 1 and 2.

\medskip
\underline{\emph{Proof of Proposition 4.} }

\medskip
\noindent Since  $V_{\tilde{\psi_C}}(i)=-i(2C-1) \approx -i 0,83193$,  \
 (\mbox{ $C$ is Catalan  constant $\approx$ 0.9159}),
then  in (1), $a_{\psi_C}=0$,  and for a measure $m_{\psi_C}$ we have  $m_{\psi_C}(\Rset)= 2C-1$. Using  (4) and  the integral representation  for $\zeta(2, s)$  function in Appendix (C),
we have
\begin{multline*}
\mathfrak{L}[ \hat{m_{\psi_C}}(x) - (2C-1)\cosh(x) ;t]= \frac{t^2/2(\zeta(2,t/2)-1/2 \zeta(2,t/4)) +1 }{t^2-1} \\ = \frac{1}{t^2-1} + \frac{1}{2} \frac{t^2}{t^2-1} \Big[ 8\int_0^\infty\frac{w}{1-e^{-2w}}e^{-tw}dw -16 \int_0^\infty\frac{w}{1-e^{-4x}}e^{-t w} dw \Big] \\  =  \mathfrak{L}[\sinh x ; t]  + 2(1+\frac{1}{t^2-1})\int_0^\infty w\frac{e^{-2w}-1}{1-e^{-4w}}e^{-t w}dw= \\ \mathfrak{L}[\sinh x ; t]+ (1+\frac{1}{t^2-1})\int_0^\infty w\frac{1- e^{2w}}{\sinh(2w)} \,e^{-t w}dw \\ =\mathfrak{L}[\sinh x  + x\frac{1- e^{2x}}{\sinh(2x)}  ; t] +\mathfrak{L}[\sinh x;t]\,  \mathfrak{L}[x \frac{1-e^{2x}}{\sinh (2x)}]\\=
\mathfrak{L}[\sinh x + x \frac{1- e^{2x}}{\sinh(2x)} +h_{\hat{\psi_C}}(x) ; t].\qquad \qquad 
\end{multline*}
 where
\begin{equation}
h_{\hat{\psi_C}}(x):= (\sinh u  \ast  u \frac{1-e^{2u}}{\sinh (2u)})(x)= \int_0^x \sinh(x-u) u \frac{1-e^{2u}}{\sinh (2u)}du.
\end{equation}
Consequently from the above calculation  we get
\begin{equation}
 \hat{m_{\psi_C}}(x) = (2C-1)\cosh x +\sinh x +\frac{x(1-e^{2x})}{\sinh(2x)} + h_{\hat{\psi_C}}(x),  \ \mbox{for $x\ge 0$.}
\end{equation}
Using website www.Wolframalpha.com or just by  elementary computations, knowing that $d/dxLi_{2}(\pm 
ix) = - x^{-1}\log(1\pm ix)$, (cf. Appendix (C)),   one checks that
\begin{multline*}
\int \sinh(x-s) s \frac{1-e^{2s}}{\sinh(2s)}ds= e^{-x}(s-1)e^s \\
- i \cosh(x) \big( - Li_2(-i e^s)+Li_2(i e^s)+ s \log(1-ie^s)- s\log(1+ie^s) \big)+ const.
\end{multline*}
Since    $ \lim_{x\to 0^+}(Ei (\pm x)-\log x)= \gamma$ (Euler constant), from (13), for $x>0$,
\begin{multline*}
h_{\hat{\psi_C}}(x)= e^{-x} (s-1)e^s  \qquad \qquad  \qquad \qquad  \qquad \\
- i \cosh(x) \big( - Li_2(-i e^s)+Li_2(i e^s)+ s \log(1-ie^s)- s\log(1+ie^s) \big)|^{s=x}_{s=0}
\end{multline*}
\begin{multline*}
= e^{-x} (x-1)e^x 
- i \cosh(x) \big( - Li_2(-i e^x)+Li_2(i e^x)+ x \log(1-ie^x) \\ - x\log(1+ie^x) \big)  - [  e^{-x} (-1) 
- i \cosh(x) \big( - Li_2(-i )+Li_2(i ) ] \\ = (x - 1) - i\cosh(x) \big ( - Li_2(-i e^x)+Li_2(i e^x) 
 +x \log \frac{1-ie^x}{1+ie^x} \big) +e^{-x}  - 2C \cosh(x)\\=
 - 1+x +e^{-x} - 2 C\cosh x  - i\cosh(x) \big ( - Li_2(-i e^x)+Li_2(i e^x) -2ix \arctan(e^x)\big),
\end{multline*}
where in the last line we  use 
$ \log(1-ie^x)-\log(1+ie^x)=-2i \arctan(e^x)$; cf. Appendix C.
And finally from  (14)  we arrive at
\begin{multline*}
\hat{m_{\psi_C}}(x) = (2C-1)\cosh x +\sinh x +\frac{x(1-e^{2x})}{\sinh(2x)}  - 1+x +e^{-x} - 2 C\cosh x 
\\ - i\cosh(x) \big ( - Li_2(-i e^x)+Li_2(i e^x) -2ix \arctan(e^x) \big) \\  =
-1 +x +  \frac{x(1-e^{2x})}{\sinh(2x)}  - i\cosh(x) \big ( - Li_2(-i e^x)+Li_2(i e^x)  -2i x \arctan(e^x)\,\big)  \qquad \qquad 
\end{multline*}
 and since $(1+ (1-e^{2x})/\sinh(2x))= - \tanh(x)$  therfore this   completes a  first part of Proposition 4.

For the second one, let us  recall that  from Jurek and Yor (2004), Corollary 1 and a formula (7) on  p. 186,  that   Khintchine (finite) measure corresponding to BDCF $\psi_C$ is given by 
\[
m_{\psi{_C}}(dx)=\frac{ \pi}{4}\, \frac{x^2}{1+x^2}\,\frac{\cosh (\pi x/2)}{\sinh^2(\pi x/2)} dx, \ \mbox{on} \ \Rset. 
\]
Hence we get that 
\begin{multline*}
\hat{m}_{\psi_{\tilde{C}}}(t)=\frac{ \pi}{2} \, \int_0^\infty \cos(tx)\frac{x^2}{1+x^2} \,\frac{\cosh (\pi x/2)}{\sinh^2(\pi x/2)}dx = 
-1 + t  +  \frac{t(1-e^{2t})}{\sinh(2t)} 
 \\ - i\cosh(t) \big ( - Li_2(-i e^t)+Li_2(i e^t)  +t \log\frac{1-ie^t}{1+ie^t}\,\big), \ \mbox{for} \ t\ge 0. \qquad \qquad 
\end{multline*}
which completes a proof of Proposition 4.

\begin{rem}
\emph{  As a consequence of Proposition 4 we have a formula
\begin{multline*}
\frac{ \pi}{2}\, \int_0^\infty \cos(tx)\frac{x^2}{1+x^2}\,\frac{\cosh (\pi x/2)}{\sinh^2(\pi x/2)}dx = 
-1 + t  +  \frac{t(1-e^{2t})}{\sinh(2t)}  \\ - i\cosh(t) \big ( - Li_2(-i e^t)+Li_2(i e^t)  -2i t \arctan(e^t)\,\big), \ \mbox{for} \ t \in \Rset,  \qquad \qquad 
\end{multline*}
that might be of some interest as well.}
\end{rem}
\underline{Proof of Corollary 1.}

\medskip
In general, if $M(dx)=h(x)dx, h(x)>0$,  is a spectral measure of a selfdecomposable characteristic function, say $\phi$, then $N(dx):= - (xh(x)^\prime dx $ is a spectral measure of a characteristic function $\psi(t):=\exp[ t \big(\log\phi(t)\big)^{\prime} ]$;  so called background driving characteristic function; cf. Jurek (1997),  Corollary 1.1., p.97  or Jurek and Yor (2004), p. 183, formulae (d) and (e).

Consequently, on the level of  Khintichine (finite) measures have;
\begin{multline*}
n(dx):=\frac{x^2}{1+x^2}N(dx) \\ = - \frac{x^2}{1+x^2}(h(x)+xh^\prime(x))dx= - m(dx) - \frac{x^2}{1+x^2} x h^\prime(x)dx,
\end{multline*}
which gives a proof of Corollary 1, when applied to $h(x):= k_C(x)$.
\begin{rem}
\emph{Since a  hyperbolic sine and a hyperbolic tangent are selfdecopmosable as well we may have statements about $\tilde{S}$ and $\tilde{T}$ analogous to the one in Corollary 1, a for  hyperbolic cosine function.}
\end{rem}

\medskip
\textbf{5. Appendix.}

For a convenience  of reading we collect here some useful  facts. Boldface numbers refer to formulae from Gradshteyn and Ryzhik (1994).

\medskip
(\textbf{A})  
$ \beta(x):=\frac{1}{2}[\ \psi(\frac{x+1}{2})-\psi(\frac{x}{2})\ ],\ \beta(x) =\sum_{k=0}^{\infty}\frac{(-1)^k}{x+k}, \ -x \notin \mathbb{N}, \   \textbf{8.732(1)}.  
$

\medskip
(\textbf{B})
For the exponential integral function $Ei$ we have :
\begin{multline*}
 \ Ei(x)= - \int_{-x}^\infty \frac{e^{-t}}{t}\,dt, \  \mbox{for} \ \ x<0; \ \    \ \ \   \textbf{8.211}(1);  \\ 
\ Ei(x)= - \lim_{\epsilon \to 0}\Big [ \int_{- x}^{-\epsilon}\frac{e^{-t}}{t}dt + \int_{\epsilon}^\infty \frac{e^{-t}}{t}dt,  \ \mbox{for}\  \ x>0;   \ \ \ \textbf{8.211}(2) \qquad \\
\ \ Ei(x)= \gamma + \log x+ \sum_{k=1}^\infty\frac{x^k}{k \cdot k!}, \ \mbox{for }x>0;  \ \ \textbf{8.214}(2); \  \  (\gamma \  \mbox{is  Euler constant})
\end{multline*}
 From above  we get: 
$
  \ \frac{d}{dx}Ei(\pm x)=\frac{e^{\pm x}}{x}$ and  $ \lim_{x\to 0^+}\,Ei (\pm x)-\log x)= \gamma.
$

\medskip
\medskip
(\textbf{C})
For  Riemann $\zeta$ functions  we have:
\[
 \zeta(s,a):= \sum_{k=0}^\infty \frac{1}{(k+a)^s}, \ \Re s > 1,\ - a \not\in \mathbb{N};  \ \ \textbf{9.521}(1) \,
\]
In particular, for $s=2$, we have $\zeta$ function as a Laplace transform:
\[
 \zeta(2,a) = \int_0^\infty \frac{x}{1-e^{-x}}e^{-a x} dx, \ \mbox{for} \ \Re a>0;  \ 
\]
\[
\mbox{because} \ \  \sum_{k=0}^\infty 
x e^{-kx}=\frac{x}{1-e^{-x}};  \ \ 
\mbox{and}  \  \mathfrak{L}[xe^{-kx}, w] = \frac{1}{(w+k)^2}.
\]
Polylogarithmic functions $Li_n(z)$, ($n$ is a fixed parameter) are defined  via series 
$$
 Li_n(z):=\sum_{k=1}^\infty \frac{z^k}{k^n}\equiv z \Phi(z,n,1) \  z \in \Cset;   \ \ \mbox{(Lerch function)}, \textbf{9.550}
$$
In particular,
$
\frac{d}{dz}Li_2(z)= -z^{-1} \log(1-z); 
  \ \ Li_2(i)-Li_2(-i)= 2 i C,
$
where  $C$ stands for a  Catalan constant.

\medskip
(\textbf{D})  For the equality, \ \ 
$
i \log \frac{1-ie^x}{1+ie^x}=2 \arctan(e^x),  \ x \in \Rset, 
$
 firstly,  note that for $x=0$,  indeed  we have   $ i \log \frac{1-i}{1+i}=\pi/2 $, and secondly,  note
equality of derivatives $ \frac{d}{dx}(i \log \frac{1-ie^x}{1+ie^x}) = \frac{d}{dx}(2 \arctan(e^x))$.

\medskip
\medskip
REFERENCES.

\medskip
[1] N.I. Akhiezer (1965),\emph{The classical moment problem}, Oliver $\&$ Boyd, 

Edinburgh and London,

\medskip
[2] L. Bondesson (1992), \emph{Generalized gamma convolutions and related}

\emph{ classes of distributions and densities},
Springer-Verlag, New York; Lecture 

Notes in Statistics, vol. 76.

\medskip
 [3] I.S. Gradshteyn and I. M. Ryzhik (1994), \emph{ Tables of Integrals, Series,}

\emph{ and Products}, Academic Press, New York.

\medskip
[4]  L. Jankowski and Z. J. Jurek (2012), Remarsk on restricted Nevalinna 

transforms, \emph{Demonstratio. Math.}, vol. XLV, no.2,
pp.297-307.

\medskip
[5]  Z. J. Jurek (1993), \emph{Operator-limit distributions in probability theory},

 J. Wiley, New York.

\medskip

[6] Z. J. Jurek (1997), Selfdecomposability:an exception or a rule ?, \emph{Ann.}

\emph{ Unviversitatis M. Curie-Sklodowska, Lublin-Polonia}, vol. LI. 1 Sectio A.

\medskip
[7 ] Z. J.  Jurek (2019), On a relation between classical and free infinitely 

divisible transforms, \emph{Probab.Math. Statist.}, to appear; 

[Also: math.arXiv:1707.02540 [math.PR],9 July 2017.]

\medskip
[8]  Z.J. Jurek and W. Vervaat (1983), An integral representation for 

selfdecomposable Banch space valued random variables.,
\emph{ Z. Wahrschein}

\emph{lichkeitstheorie verv. Gebiete}, vol. 62, pp. 247-262.

\medskip
[9] Z. J. Jurek and M. Yor (2004), Selfdecomposable laws associated with 

hyperbolic functions,\emph{Probab. Math. Statist.}, vol.24, Fasc.1, pp. 181-190.

\medskip
[10] S. Lang(1975), \emph{$SL_2(\Rset)$}, Addis0n-Wesley, Reading Massachusetts.

\medskip
[11]  J. Pitman and M. Yor (2003),  Infinitely divisible laws associated with 

hyperbolic functions,\emph{ Canad. J.
Math. } 55 (2), pp. 292-330.

\medskip
\medskip
Author's address: 

Institute of Mathematics, University of Wroc\l aw, Pl. Grunwaldzki 2/4, 

50-384 Wroc\l aw, Poland;

Email:  zjjurek@math.uni.wroc.pl ;   www.math.uni.wroc.pl/$\sim$zjjurek

\end{document}